\magnification=\magstep1
\input amstex
\UseAMSsymbols  
\hoffset=0truecm \voffset=0truecm 
\NoBlackBoxes
   
\font\rmk=cmr8    
\input pictex

\def\mat{\operatorname{mat}}

\def\Hom{\operatorname{Hom}}

\def\add{\operatorname{add}}

\def\Ext{\operatorname{Ext}}

\def\mod{\operatorname{mod}}
\def\End{\operatorname{End}}

\def\arr#1#2{\arrow <1.5mm> [0.25,0.75] from #1 to #2}
\def\op{{\text{op}}}
\vglue1cm

\centerline{\bf Dieter Vossieck and the Development}
\centerline{\bf of the Representation Theory of Artin Algebras}
	\medskip
\centerline{Claus Michael Ringel}
	\bigskip\bigskip

This is a pre-dinner-lecture, as part of the Bielefeld workshop 
{\it Discrete Categories in Representation Theory}, April 20 -- 21, 2018.
The workshop discusses the possible discreteness of additive
categories: starting point is the 
celebrated paper {\it The algebras with discrete derived category} [9]
by Dieter Vossieck. But there is another (more secrete) aim for this workshop:
to celebrate the 60th birthday of Dieter Vossieck --- he was
born in September 1957. In this lecture, I want to give a short report 
on some of his contributions to representation theory, concentrating on his early
Bielefeld years 1978 -- 1986. Dieter Vossieck may be considered as one of my best
PhD students, better: he is it, but unfortunately, he never handed in what was 
supposed to be his Bielefeld PhD thesis (see the comments below). 
Dieter provided a lot of important contributions to representation theory ---
actually there are only few publications, but his influence is much wider. The aim of my
lecture is to draw the attention not only to his publications, but also to some of his
ideas which he presented in lectures and in discussions and which are not available
otherwise. He always was a perfectionist, so he refused to put forward any incomplete 
or unpolished result. 

Dieter Vossieck studied at Bielefeld University, starting in 1978, he stayed here for
eight years, then he went to Z\"urich, this was in 1986. 
He worked in Switzerland (in Z\"urich as well as in Basel) again 
for eight years. Then, from 1994 to 2008 he was in Mexico, at the UNAM in Mexico City
and at the Universidad Michoacana des San Nicol\'as de Hidalgo in Morelia,
interrupted for one year (2001 -- 2002) by being professor at Beijing Normal University,
in China.
Since 2008 he is back in Bielefeld, partly as a private scholar, partly with some temporary 
university contracts.

	\bigskip\bigskip
{\bf 1\. The Happel-Vossieck list.}  
	\medskip
The first topic to be mentioned concerns the joint paper of Happel and
Vossieck  entitled {\it Minimal algebras of infinite
representation type with preprojective component} [1]; 
it presents what is now the famous Happel-Vossieck list, namely 
the list of the frames of the tame concealed algebras.

Throughout the lecture, I will assume that we work over an algebraically
closed field $k$. A {\it tame concealed} algebra is the endomorphism ring 
$A = \End({}_HT)$ of a
preprojective tilting $H$-module $T$, where $H$ is a (connected) tame hereditary algebra.
It turns out that 
if $A$ is tame concealed, then also the opposite algebra $A^\op$ is tame concealed. 

Let us assume that $H$ is a tame hereditary algebra and $T$ a preprojective tilting
$H$-module with endomorphism ring $A$. 

If the quiver of $H$ is a cycle with $p+q$ arrows, where $p$ arrows point in one direction
and $q\le p$ in the opposite direction, then $H$ will be said to be of {\it type}
$\widetilde{\Bbb A}_{pq}$. In this case, $T$ has to be a slice module, thus $A$ is
hereditary again, and of the same type $\widetilde{\Bbb A}_{pq}$.
If the quiver of $H$ is not a cycle, then its underlying 
graph is one of the extended Dynkin diagrams
$\widetilde{\Bbb D}_n,\ \widetilde{\Bbb E}_6,\ \widetilde{\Bbb E}_7,\  
\widetilde{\Bbb E}_8$ and this is called the {\it type} of $A$.
If $H$ is of type $\widetilde{\Bbb D}_n$, then, up to change of orientation of
arrows, there are at most four possibilities for the quiver of $A$.
The most interesting algebras are those of type 
$\widetilde{\Bbb E}_6,\ \widetilde{\Bbb E}_7,\  
\widetilde{\Bbb E}_8$. There are 56, 437, and 3\,801 
different isomorphism classes of tame concealed algebras of type $\widetilde{\Bbb E}_6$
$\widetilde{\Bbb E}_7$, and $\widetilde{\Bbb E}_8$, respectively.
Now 3\,801 different classes may seem to be difficult to overlook, but it is possible to
reduce the numbers considerably by looking at {\it frames}. For 
$\widetilde{\Bbb E}_8$, there are 117 possible frames, 29 for $\widetilde{\Bbb E}_7$
and 5 for $\widetilde{\Bbb E}_6$,  and it is easy to obtain the
corresponding isomorphism classes up to duality by just 
replacing the arms by corresponding branches. 

What is the relevance of this classification? The tame concealed algebras
are important minimal representation-infinite algebras, where {\it minimal
representation-infinite} means that the algebra $A$ itself is representation-infinite,
whereas for any non-zero idempotent $e$, the factor algebra $A/\langle e\rangle$
is representation-finite. What Happel and Vossieck show is that
{\it an algebra $A$ with a preprojective 
Auslander-Reiten component is minimal representation-infinite if and only if
it is either tame concealed or else a wild generalized Kronecker algebra} 
(a generalized Kronecker algebra is the
path algebra of a directed quiver with two vertices and say $n$ arrows; it is wild
provided $n\ge 3$). 

There are numerous applications: First of all,  
the Happel-Vossieck list is used by Bongartz in order to
characterize the representation-finite algebras (for this reason, the list
is sometimes called the BHV-list). Secondly, the famous multiplicative basis
paper by Bautista-Gabriel-Roiter-Salmeron strongly relies on the list (thus the paper
reprints the list, but using a different ordering of the frames), as do the various 
proofs of the second Brauer-Thrall conjecture. But the list plays a role also
in other settings, for example in the context of cluster categories. 

In 1984, Peter Gabriel gave a colloquium lecture at Bielefeld, devoted to the foundation
of the representation theory of finite-dimensional algebras. At the end, he suggested that 
the Happel-Vossieck list has a similar relevance for representation theory
as the list of the finite simple groups has for group theory. I should recall that 
at that time the classification of the finite simple groups was very popular --- it
was considered as 
one of the main achievements of the century; it was discussed even in newspapers. 
Gabriel's rating was felt as a real provocation at a mathematical faculty which was nicknamed
outside of Bielefeld as a
{\it Faculty of Group Theory} (since nearly all major mathematicians at 
Bielefeld were working on 
problems related to group theory, and at least one of them, Bernd Fischer, was
strongly involved in the classification project). 
	\medskip 
The Happel-Vossieck 
classification was fully presented at a Luminy conference in September 1982 
(and not {\it thereafter}, as formulated by Bongartz in {\it manuscripta mathematica}
vol 46 (1984); actually there was a big discussion comparing the 
Bongartz and the Happel-Vossieck approach). It was published in 1983 and it is the 
center of Vossieck's Diplom thesis at Bielefeld. The
thesis was handed in only in 1984, since Dieter 
was unhappy with a proof which relies on the use of a computer. Also,
being a perfectionist, Dieter felt that the presentation of the frames
should be based on some intrinsic partial ordering. 
He often came to my office showing me huge
pieces of papers visualizing his attempts to outline such an ordering. 
Indeed, already a first glance at the various frames confirms the idea of a partial ordering of
the tame concealed algebras of a fixed type, with the hereditary algebras as the maximal elements
and the unique canonical algebra as the minimal element. 
Given a tame concealed algebra $A$, the Euler
form $q_A$ on the Grothendieck group of the module category is semidefinite with a positive radical
generator $\bold h$ and the Happel-Vossieck list provides for each frame such a generator
$\bold h$. The sum of the coefficients of $\bold h$ is an interesting invariant: it is maximal
for $A$ hereditary and minimal for the canonical algebra $A$.
A similar invariant is what I want to call the {\it shortage} of $A$, it is the number of
roots of $q_A$ which are neither positive, nor negative: the shortage of $A$ is zero iff
$A$ is hereditary, and is maximal iff $A$ is the canonical algebra. 
Until now, no definite ordering of the tame concealed algebras has been described in detail;
that is really a desideratum.  
	\bigskip\bigskip
{\bf 2. The foundation of tilting theory.}
	\medskip
It seems that the first presentation of the Brenner-Butler tilting theory was given at
a two day workshop in Bielefeld, in spring 1979, by Butler. His lecture was 
on the first day, and after this lecture, there was the immediate wish to see more
details, thus a continuation was scheduled for the second day.
Of course, special tilting functors were known at that time (first the Bernstein-Gelfand-Ponomarev
reflection functors, then the APR-tilts, but also other 
constructions which were used for example 
by students of Gabriel). What was new and exciting was the axiomatic approach exhibited by
Butler. It soon turned out that these axioms could be weakened considerably without 
changing the essential results (the Brenner-Butler tilting modules $T$ had a projective
direct summand which generates $T$). 

What is called the Brenner-Butler theorem [2] is usually illustrated in the  following picture.
Here, we start with an $A$-module $T$ which is {\it tilting}: it has projective dimension at most 1,
no self-extension, and there is an exact sequence $0 \to {}_AA \to T' \to T'' \to 0$
with $T',T''$ in $\add T$. Let $B = \End(T)^\op$.
$$
\hbox{\beginpicture
\setcoordinatesystem units <.9cm,1cm>
\put{} at 1 0
\put{} at 13.7 6

\plot 2.5 4.5  2.7 4.7  2.5 4.9  2.5 5.3  2.7 5.5 /
\plot 2.5 4.5  4.5 4.5  5 4  7 4  7.5 4.5  9.5 4.5 / 
\plot 9.5 4.5  9.7 4.7  9.5 4.9  9.5 5.3  9.7 5.5 /
\plot 2.7 5.5  4.5 5.5  5 6  7 6  7.5 5.5  9.7 5.5 /

\plot 9.5 0.5  9.7 0.7  9.5 0.9  9.5 1.3  9.7 1.5 /
\plot 9.7 0.5  9.9 0.7  9.7 0.9  9.7 1.3  9.9 1.5 /
\plot 6 2  7 2  7.5 1.5  9.7 1.5 /
\plot 6.2 0  7 0  7.5 0.5  9.5 0.5 /

\put{$*$} at 6 6
\put{$*$} at 6.5 5.5
\put{$*$} at 5 5
\put{$*$} at 6.7 4.5
\put{$*$} at 6.2 4

\setquadratic
\plot 6 6      6.65  5.75     6.5 5.5  /
\plot 6.5 5.5  6.55  5.25    5 5  /
\plot 5 5      6.55  4.75   6.7 4.5 /
\plot 6.7 4.5  6.85  4.25   6.2 4 /

\plot 5.8 6    5.65 5.75     6.3 5.5  /
\plot 6.3 5.5  4.85 5.25  4.8 5  /
\plot 4.8 5    5.05 4.75    6.5 4.5 /
\plot 6.5 4.5  6.05 4.25   6 4 /

\plot 6 2      6.65  1.75     6.5 1.5  /
\plot 6.5 1.5  6.55  1.25    5 1  /
\plot 5 1    6.55  0.75   6.7 0.5 /
\plot 6.7 0.5  6.85  0.25   6.2 0 /

\plot 13 2    12.85 1.75     13.5 1.5  /
\plot 13.5 1.5  12.05 1.25  12 1  /
\plot 12 1    12.25 0.75    13.7 0.5 /
\plot 13.7 0.5  13.25 0.25   13.2 0 /

\setlinear
\plot 9.9 1.5  11.7 1.5  12.2 2  13   2 / 
\plot 9.7 0.5  11.7 0.5  12.2 0  13.2 0 /
\put{$\mod A$} at 1 5
\put{$\mod B$} at 1 1

\put{$\Cal F$} at 4 5
\put{$\Cal T$} at 8 5
\put{$\Cal Y$} at 8 1
\put{$\Cal X$} at 11 1


\arr{8 4}{8 2}
\arr{10.98 2.1}{11 2}
\put{$\Hom_A(T,-)$} at 9.15  3.6 
\put{$\Ext^1_A(T,-)$} at 3.1 3.6

\setquadratic
\setdashes <1.5mm>
\plot 4 4  4.4 3.3  5 3 7.5 3  7.7 3 /

\plot 8.1 3  9 3  10 3  10.6 2.7  11 2 /

\endpicture}
$$
The class $\Cal T$ is the class of $A$-modules
generated by $T$ and the class $\Cal Y$ is the class of $B$-modules cogenerated by 
the dual $D(T)$ of $T$. What Brenner and Butler had observed, and this is 
the first step of tilting theory, is that the functor $\Hom_A(T,-)$ provides
an equivalence between $\Cal T$ and $\Cal Y$ (the vertical arrow in the center of the
picture). 

But the tilting module $T$ provides also a second equivalence, namely the equivalence
(given by $\Ext^1_A(T,-)$ and shown by a bended dashed arrow) 
between the class $\Cal F$ of all $A$-modules $M$ with
$\Hom(T,M) = 0$ and the class $\Cal X$ of all $B$-modules $N$ with
$\Hom(N,D(T)) = 0$. After Butler's Bielefeld lecture, a lot of examples were
studied in Bielefeld, and it took quite a while to realize that one always has 
the second equivalence. 
As I remember, the examples calculated by Dieter Vossieck were those
which really were convincing and which paved the way to the final proof. 

In this way, Dieter Vossieck was involved right from the beginning in the development
of tilting theory. In the proceedings of ICRA III, Klaus Bongartz published an
account of tilting theory which ended with the (slightly presumptuous) remark
that {\it well-read mathematicians tend to understand the tilting theorem using spectral
sequences.} Of course, this concerned the classical tilting modules of projective
dimension at most 1. Thus, when the general notion of tilting modules with arbitrary
finite projective dimension was introduced and studied by Miyashita, Happel, as well as  Cline-Parshall-Scott, Dieter Vossieck (who always was well-read) 
took the initiative to analyze the corresponding tilting functors using 
spectral sequences. This was supposed to be his PhD-thesis at Bielefeld. In the summer 1986
he gave several lectures on his work. The manuscript was completely finished before
he went to Z\"urich. Unfortunately, he hesitated to hand it in ---
apparently he wanted to show it first to Gabriel (but already in 1984, the relationship
between Z\"urich and Bielefeld had started to be quite frosty). At the time, when Dieter 
was leaving Z\"urich in order to go to Mexico, I was invited to a lecture at Z\"urich. As usual,
after such a lecture there was a common dinner and during the dinner Gabriel stressed that 
he did not understand why Dieter had not handed in the Bielefeld thesis
(I have to say, he said this to my surprise, but may-be also Dieter's).

More than 15 years later, 
Brenner and Butler provided a contribution [11] for the {\it Handbook of Tilting Theory}
with the title {\it A spectral analysis of classical tilting functors.} There, they
write: {\it The spectral sequences in question seem first to have been written
down, but nowhere published, by Dieter Vossieck in the mid-1980's and
were re-discovered by the authors during the summer of 2002 whilest preparing the
talk for the conference ``Twenty Years of Tilting Theory'' at Chiemsee in November 2002
on which this article is based. After that talk, Helmut Lenzing mentioned Vossieck's work, and
kindly supplied a copy of his notes of a lecture in July 1986 by Vossieck at the 
University of Paderborn entitled ``Tilting theory, derived categories and spectral
sequences''} [... ;] {\it in the
last two sections Vossieck briefly described the spectral sequences and
filtration formulae which are stated and proved in the main part, Section 3, of this article.}
In this way, a spectral sequence approach to tilting theory was finally made available
(whereas the original manuscript of Dieter seems to be lost, as he told me ---
that is a pity). 

	\bigskip\bigskip
{\bf 3\. Further topics.}  
	\medskip
As I mentioned at the beginning, this report is devoted mainly to the early Bielefeld years
of Dieter Vossieck, but let me draw the attention at least to some of the topics
he later was working on. Of course, he continued to look at questions in tilting theory:
there are three very important papers written with Bernhard Keller, two in the
Comptes Rendus [4,5], one in the Bulletin of the Belgian Mathematical Society [6], all three
concern triangulated categories and aim to clarify and to complement the work of 
Dieter Happel, in particular his lecture notes on the use of triangulated categories in the
representation theory of finite-dimensional algebras. Since tilting theory for
module categories deals with special torsion pairs, one should, in the context of
triangulated categories, investigate t-structures. And this is, what they do.
Comparing the hearts of different
t-structures, Keller and Vossieck 
obtained a common generalization of the Grothendieck-Roos duality for 
regular commutative rings and the tilting theory for finite-dimensional
algebras. Also, tilting was generalized by Keller-Vossieck to silting (and silting theory
became really popular in recent years).
	\medskip 
Vossieck obtained his PhD at the University of Z\"urich, with a thesis [8] dealing with
matrix problems (and not related at all to his Bielefeld PhD project). At that time,
Nazarova and Roiter were very eager to be
in contact with mathematicians from the West.
They had developed many powerful, but technical 
reduction methods for matrix problems and hoped that for instance Gabriel would provide 
some fancy categorical or homological interpretation. The paper {\it Tame and wild
subspace problems} [7] by Gabriel, Nazarova, Roiter, Sergejchuk and Vossieck has to be
seen as part of this Kiev-Z\"urich cooperation.
	\medskip
For many applications, in particular for applications in Lie theory and in the theory
of algebraic groups, quasi-hereditary algebras play a decisive role. 
Given a quasi-hereditary algebra $A$, say with the set $\Delta$ of standard modules,
one is interested in the category $\Cal F(\Delta)$ of all modules with a filtration
with factors in $\Delta$.
Starting with a bimodule $B$, one may define the category of matrices $\mat B$ over $B$.
In case $B$ is what is called upper-triangular, $\mat B$ can be identified with
the category $\Cal F(\Delta)$ for some quasi-hereditary algebra $A$.
Now the radical of a hereditary finite-dimensional
algebra is upper-triangular, thus $\mat B = \Cal F(\Delta)$ for some quasi-hereditary algebra 
$A$. A joint publication [10] with Hille provides in this case an explicit description of the
quasi-hereditary algebra $A$. 
	\medskip
Of course, there is the 2001 paper [9] on the algebras with discrete derived categories, which
has attracted a lot of interest in the last years and which is the basis for the present
workshop. Vossieck classified these algebra. Examples are the piecewise
hereditary algebras of Dynkin type. Up to Morita equivalence, the remaining 
algebras with discrete derived categories 
are gentle algebras with a unique cycle in the quiver, and with an additional condition
on the relations (the so-called clock condition). Soon after (in 2004) 
Bobi\'nski, Gei\ss{} and Skowro\'nski presented representatives for the derived
equivalence classes of these algebras: besides the Dynkin algebras, there are the 
gentle algebras $\Lambda(r,n,m)$, with $1\le r \le n$ and $0 \le m$; they are
given by an oriented cycle with $n$ arrows and $r$ zero relations of length 2,
with an arm of length $m+1$ attached. They also described the Auslander-Reiten
quiver of the derived category: there are $2r$ components of type $\Bbb Z \Bbb A_\infty$
and $r$ components of type $\Bbb Z \Bbb A_\infty^\infty$. These algebras and categories
have been studied further by several mathematicians and now are quite well
understood. 
	\medskip
In 2009, Dieter Vossieck organized at Bielefeld a reading course 
on ray categories (in particular about the freeness of the fundamental group, 
Roiter's vanishing theorem, and interval-finiteness) and this was followed by
a sequence of lectures by him on multiplicative bases for 
subspace-finite vector space categories. Ray categories as introduced
by Bautista, Gabriel, Roiter and Salmeron in their multiplicative basis paper
have to be seen as the essential language for dealing with finite-dimensional
algebras whose module categories are not too complicated or not too much entangled
(whatever this means). The multiplicative basis paper and some consecutive papers
provide a lot of fundamental structure theorems. 
Unfortunately, these papers are not so easy to follow, since they are very technical.
They have scared away a lot of people (and the many recent books devoted to the
representation theory of finite-dimensional algebras just avoid to touch these questions).
Since no comprehensive presentation is yet available, 
any publication which provides help, provides improvement would be strongly appreciated.
When Dieter delivered the lectures, he promised a written
account, with all his fascinating examples which illuminate the
results as well as the difficulties. But we are still waiting \dots. 
Such an account, even if it is
incomplete and covers only specific parts, would be really valuable and definitely
appreciated by the mathematical community.
	\bigskip\bigskip
{\bf 4. Hammocks}
	\medskip
There is just one joint publication of Dieter and myself. I like it very much, 
but it seems to be nearly forgotten --- thus I want to use this
occasion as an advertisement. The paper is called {\it Hammocks} [3] and has been 
published in the Proceedings of the London Mathematical Society in 1987.
The concept of a hammock had been introduced before by Sheila Brenner in order to provide
a combinatorial characterization of the function which attaches to any indecomposable
module $M$ the Jordan-H\"older multiplicity $[M\colon E]$ of some fixed simple module $E$. 
	\medskip
A {\it hammock} is a finite translation quiver $H$ with a unique source $\omega$
and with an additive function $h$ with positive values
(the {\it hammock function}) such that
$h(x) = 1$ for any vertex $x$ which is projective or injective. 
A typical example is the following: Let $A$ be a representation-directed algebra
and $E$ a simple module. Then the full subquiver of the Auslander-Reiten quiver 
$\Gamma(A)$ of $A$
given by the indecomposable modules $M$ with $[M\colon E]\neq 0$ is a hammock with hammock
function $h(M) = [M\colon E]$ and with $\omega$ being the projective cover of $E$ (and $\Gamma(A)$ 
itself is just the union of these hammocks).   
	\medskip
The main result of our hammock paper is as follows: {\it There is a bijection between
the hammocks $H$ and the 
subspace-finite posets $X$, where $H$ corresponds to the
poset $X$ of all projective vertices $p\in H$ different from $\omega$. 
Conversely, one attaches to a subspace-finite poset $X$ the Auslander-Reiten quiver $H$ 
of its subspace category.} 
	\medskip
This shows that any hammock arising in whatever part of mathematics provides an
interpretation of the corresponding setting in terms of subspace categories: a purely
combinatorial invariant, the hammock function, points to a realization using subspaces of
vector spaces. 
	\medskip
Let us formulate an immediate consequence. 
We call an $A$-module $M$ {\it locally indecomposable} provided for any primitive
idempotent $e$ of $A$, 
the vector space $eM$ with its definable subspaces is an indecomposable $S$-space.
(Recall that an {$S$-space} $(V;U_i)_{i\in I}$ is given by a vector space $V$ and a set
of subspaces $U_i$ of $V$ which are indexed by some set $I$. Such an $S$-space $(V;U_i)_{i\in I}$
is said to be {\it indecomposable} provided there is no proper decomposition 
$V = V'\oplus V''$ with
$(U_i\cap V')+(U_i\cap V'') = U_i$, for all $i\in I$.) 
	\medskip 
{\bf Corollary.} {\it Let $A$ be a representation-directed algebra. Then
any indecomposable $A$-module $M$ is locally indecomposable.}
	\medskip
Let us consider some examples of indecomposable modules $M$ and primitive idempotents $e$ such that
$eM$ is 2-dimensional. To say that $eM$ together with some subspaces $U_j$ is
an indecomposable $S$-space means in this case that (at least) 3 of these subspaces $U_j$ are
1-dimensional and pairwise different. 

Here are the examples: we draw quivers with relations (a commutative relation in the left
column, a zero relation in the right one) and exhibit a dimension vector $\bold d$. The
module to be considered is 
the (uniquely determined) indecomposable representation $M$ with dimension vector $\bold d$. 

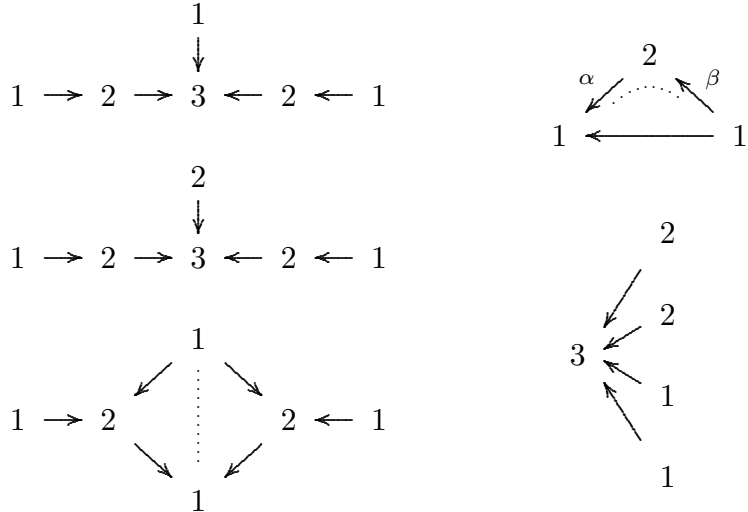
$$
{\beginpicture
\setcoordinatesystem units <1cm,.9cm>
\put{\beginpicture
\multiput{$1$} at 0 0  4 0  2 1 /
\multiput{$2$} at 1 0  3 0  /
\multiput{$3$} at 2 0 /
\arr{0.3 0}{0.7 0}
\arr{1.3 0}{1.7 0}
\arr{2.7 0}{2.3 0}
\arr{3.7 0}{3.3 0}
\arr{2 0.7}{2 0.3}
\endpicture} at 0 1
\put{\beginpicture
\multiput{$1$} at 0 0  4 0 /
\multiput{$2$} at 1 0  3 0  2 1 /
\multiput{$3$} at 2 0 /
\arr{0.3 0}{0.7 0}
\arr{1.3 0}{1.7 0}
\arr{2.7 0}{2.3 0}
\arr{3.7 0}{3.3 0}
\arr{2 0.7}{2 0.3}
\endpicture} at 0 -1
\put{\beginpicture
\multiput{$1$} at 0 0  4 0  2 1  2 -1 /
\multiput{$2$} at 1 0  3 0 /
\arr{0.3 0}{0.7 0}
\arr{1.7 0.7}{1.3 0.3}
\arr{3.7 0}{3.3 0}
\arr{1.3 -.3}{1.7 -.7}
\arr{2.7 -.3}{2.3 -.7}
\arr{2.3 0.7}{2.7 0.3}
\setdots <1mm>
\plot 2 0.6  2 -.6 /
\endpicture} at 0 -3.5
\put{\beginpicture
\multiput{$1$} at 0 0  2 0 /
\put{$2$} at 1 1
\arr{0.7 0.7}{0.3 0.3}
\arr{1.7 0.3}{1.3 0.7}
\arr{1.7 0}{0.3 0}
\setdots <1mm>
\setquadratic
\plot 0.6 0.4  1 0.6  1.4 0.4 /
\put{$\ssize \alpha$} at 0.3 0.7 
\put{$\ssize \beta$} at 1.7 0.7 
\endpicture} at 5 0.5 
\put{\beginpicture
\multiput{$2$} at 1 1.5  1 0.5 /
\multiput{$1$} at 1 -.5  1 -1.5 /
\multiput{$3$} at 0 0  /
\arr{0.7 1.05}{0.3 .35}
\arr{0.7 0.35}{0.3 0.075}
\arr{0.7 -1.05}{0.3 -.35}
\arr{0.7 -0.35}{0.3 -.075}
\endpicture} at 4.7 -2.7
\endpicture}
$$
On the left, we deal with representation-directed algebras,
thus the Corollary can be applied: for any primitive idempotent $e$, the
vector space $eM$ with its definable subspaces is an indecomposable $S$-space.

In contrast, on the right,
the 2-dimensional vector spaces (with their definable subspaces)
are decomposable. In the upper case, the only non-zero proper subspace which is
definable is the image of $\beta$ (which is equal to the kernel
of $\alpha$). In the lower example, the 2-dimensional spaces again
have only one non-zero proper subspace which is definable,
namely (if we assume that we deal with subspaces and inclusion maps) the intersection
of the two 2-dimensional subspaces. 
	\medskip
In this way, one obtains an effective way to attach to a faithful indecomposable $R$-module $M$ 
of a representation-finite basic algebra $R$ its covering $\widetilde R$-module 
$\widetilde M$, where $\widetilde R$ is the universal cover of $R$.  Namely, let
$M = \bigoplus M_i$ be the Peirce decomposition of $M$ (thus $M_i = e_iM$, where 
$e_i,\dots,e_n$ is a complete set of pairwise orthogonal, 
primitive idempotents of $R$). 
Consider $M_i$ together with its definable subspaces and write it as the direct sum
of indecomposable $S$-spaces $M_{ij}$. Then $\widetilde M = \bigoplus_{ij} M_{ij}$
is the Peirce decomposition of $\widetilde M$. 

Above, we have exhibited an indecomposable representation $M$ of a
quiver with three vertices, three arrows and a zero relation
(the upper case in the second column). Let us show the corresponding 
covering module $\widetilde M$:
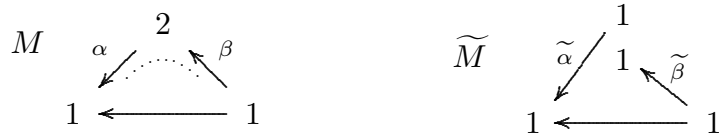
$$
{\beginpicture
\setcoordinatesystem units <1cm,1cm>
\put{\beginpicture
\multiput{$1$} at 0 0  2 0 /
\put{$2$} at 1 1
\arr{0.7 0.7}{0.3 0.3}
\arr{1.7 0.3}{1.3 0.7}
\arr{1.7 0}{0.3 0}
\setdots <1mm>
\setquadratic
\plot 0.6 0.4  1 0.6  1.4 0.4 /
\put{$\ssize \alpha$} at 0.3 0.7 
\put{$\ssize \beta$} at 1.7 0.7 
\put{$M$} at -.5 0.8
\endpicture} at 0 0
\put{\beginpicture
\multiput{$1$} at 0 0  2 0  1 0.7  1 1.2 /
\arr{0.8 1}{0.25 0.25}
\arr{1.7 0.2}{1.2 0.6}
\arr{1.7 0}{0.25 0}
\put{$\ssize \widetilde\alpha$} at 0.35 0.78 
\put{$\ssize \widetilde\beta$} at 1.6 0.6 
\put{$\widetilde M$} at -.7 0.8
\endpicture} at 5 0 
\endpicture}
$$

	\medskip
{\bf Finale.} In 1996, Bielefeld celebrated the 60th birthday of Bernd Fischer.
One of the major speakers was John Horton Conway, well-known not only for his 
contributions to group theory, but also for his interest in numbers and games.
He stayed at Bielefeld for nearly a week, always having a puzzle, the 
so-called Hanayama Cast Devil, in his pocket and 
asking anyone he was speaking to to separate the two identical 
(and quite innocent looking) metal pieces
--- a surprisingly difficult task. 
Actually, one may increase the difficulty by using three instead of the
usual two pieces. It seems to me that this may be a suitable birthday present for Dieter.
Happy Birthday!
	\bigskip
{\bf References}
	\medskip
\frenchspacing 
\item{[1]}  D. Happel, D. Vossieck. Minimal algebras of infinite representation type 
   with preprojective component. manuscripta math. 42 (1983), 221--243.
\item{[2]}  D. Happel, C. M. Ringel. Tilted algebras. Trans. Amer. Math. Soc. 274 (1982), 399--443. 
\item{[3]}  C. M. Ringel, D Vossieck. Hammocks. Proc. London Math. Soc. (3) 54 (1987), 216--246. 
\item{[4]}  B. Keller, D. Vossieck.  Sous les ca\'egories d\'eriv\'ees. C. R. Acad. Sci.
   Paris. S\'er. I. 305 (1987), 225--228.
\item{[5]}  B. Keller, D. Vossieck.  Dualit\'e de Grothendieck-Roos et basculement
   C. R. Acad. Sci. S\'er I. 307 (1988), 543–-546.
\item{[6]}  B. Keller, D. Vossieck. Aisles in derived categories.
   Bull. Soc. Math. Belg. 40 (1988), 239--253.
\item{[7]}  P. Gabriel, L. A. Nazarova, A. V. Roiter, V. V. Sergejchuk, D. Vossieck.
   Tame and wild subspace problems, Ukr. Math. J. 45 (1993), 313--352.
\item{[8]} D. Vossieck. Representation-finite weakly completed posets. Dissertation. 
   Universit\"at Z\"urich (1993).
\item{[9]}  D. Vossieck. The algebras with discrete derived category.
   Journal of Algebra 243 (2001), 168--176.
\item{[10]}  L. Hille, D. Vossieck.
   The quasi-hereditary algebra associated to the radical bimodule over a hereditary algebra
  Colloquium Mathematicae 98 (2003). 201--211.
\item{[11]}  S. Brenner, M. C. R. Butler. A spectral analysis of classical tilting functors.
  in: Handbook of Tilting Theory
  (ed. L. Angeleri H\"ugel, D. Happel, H. Krause). London Math. Soc. Lecture Note Series vol 332.
  Cambridge University Press (2007). 31--48.

	\bigskip 
{\rmk
C. M. Ringel\par
Fakult\"at f\"ur Mathematik, Universit\"at Bielefeld \par
POBox 100131, D-33501 Bielefeld, Germany  \par

\bye
\vfill\eject
At the Luminy conference in 1982, Happel distributed the Happel-Vossieck list just as it 
later has been published for example in Springer LNM 1099). Actually, it was
a challenge to squeeze the printout so that it fitted on a single page --- 
I helped him, working at a copy machine, using scissors and glue ... . At Luminy, Happel
have two lectures (of one his topics was the list) as did Bongartz. Bongartz presented
his criterion, but referred to a much larger class of algebras. When we told him (in between
his two lectures) that the
HV-list should do the job, he strongly objected and insisted again in the second lecture 
that his criterion needs the larger class of algebras and that it would be impossible to present then in a similar concise way.